\documentclass[letterpaper, 10 pt, journal, twoside]{IEEEtran}	

\IEEEoverridecommandlockouts                              

\usepackage{comment}
\usepackage{graphics} 
\usepackage{epsfig} 
\usepackage{mathptmx} 
\usepackage{times} 
\usepackage{amsmath,bm,mathtools} 
\usepackage{amssymb}  
\usepackage[utf8]{inputenc}
\usepackage[T1]{fontenc}
\usepackage[section]{placeins}
\usepackage{dsfont}
\usepackage{algpseudocode}
\usepackage{algorithm}
\usepackage{url}
\usepackage{tabularx}
\usepackage[percent]{overpic}
\usepackage{cite} 
\usepackage{calrsfs}
\DeclareMathAlphabet{\pazocal}{OMS}{zplm}{m}{n}
\usepackage{color}

\newtheorem{definition}{Definition}
\newtheorem{assumption}{Assumption}

\newtheorem{lemma}{Lemma}
\newtheorem{thm}{Theorem}

\newcommand{\col}{\mathrm{col}}

\newcommand{\R}{\mathbb{R}}

\newcommand{\bs}[1]{\boldsymbol{#1}}

\usepackage{xcolor,calc}

\newcommand{\colorNoteBox}[3]{
	\begin{center}
		\vspace{-2ex}\small
		\fcolorbox[rgb]{#1}{#2}{\parbox[t]{\columnwidth-\parindent- 0.2cm}{\setlength{\parskip}{1.5ex}#3}}
	\end{center}
}
\newcommand{\FFnote}[1]{\colorNoteBox{0,1,0}{0.9,0.9,0.9}{FF: #1}}

\begin{document}
%

\title{Probabilistically robust stabilizing allocations\\ in uncertain coalitional games}

\author{
	George Pantazis, Filippo Fabiani, Filiberto Fele and Kostas Margellos
	\thanks{The authors are with the Department of Engineering Science, University of Oxford, OX1 3PJ, United Kingdom {\tt \footnotesize \{georgios.pantazis, filippo.fabiani, filiberto.fele, kostas.margellos\}@eng.ox.ac.uk}. This work was partially supported by the UK Engineering and Physical Sciences Research Council [grant EP/P03277X/1], and Innovate UK [grant number 104781]. 
}
}

\maketitle
\thispagestyle{empty}
\pagestyle{empty}

\begin{abstract} 
	In this paper we consider multi-agent coalitional games with uncertain value functions for which we establish distribution-free guarantees on the probability of allocation stability, i.e., agents do not have incentives to defect from the grand coalition to form subcoalitions for unseen realizations of the uncertain parameter.
	In case the set of stable allocations, the so called core of the game, is empty, we propose a randomized relaxation of the core. We then show that those allocations that belong to this relaxed set can be accompanied by stability guarantees in a probably approximately correct fashion.  Finally, numerical experiments corroborate our theoretical findings. 
	
\end{abstract}


\section{Introduction}
Multi-agent systems are ubiquitous in  applications ranging from engineering \cite{Raja, Malcolm2, Fele_coalitional} to economics  and social sciences \cite{Mccain2008}.  In such systems, agents typically act as selfish entities that optimize their own payoff function. In some cases, however, due to their limited ability to increase their own utility when working on their own, agents have an incentive to form coalitions aiming at receiving a higher individual payoff. This gives rise to a coalitional game setup  \cite{Chalkiadakis}. Since each agent is interested in maximizing their own welfare, the problem of allocating the total value of the coalition in order to guarantee that none of the agents has an incentive to defect it, is key from both a collective and an individual point of view. This concept is known as \emph{stability} of the agents' allocations and it is related to the agents' coalitional values. \par
In real-world applications, the values of each coalition are typically affected by uncertainty. This can be due to various reasons. The most predominant ones refer to the effects of exogenous factors acting on the agents' network, or uncertainty inherent in the coalition formation process. This paper focuses on the former case, where  uncertainty in the environment results in changes in the values of agents' coalitions.   Encoding uncertainty in the value functions of a coalitional game was initially explored in the seminal works \cite{Charnes_Chance_I,Charnes_Chance_II,Charnes1973}. Therein, the authors focus on the extension of well-known solution concepts to provide stability-wise robust allocations. In \cite{Suijs1} it is shown that for a certain class of stochastic games, properties such as the non-emptiness of the core of a deterministic game continue to hold for their stochastic counterpart. References \cite{Chalkiadakis_RL} and \cite{Ieon_Bayesian} tackled uncertainty in the values of the coalitions by leveraging  Bayesian learning methods, while \cite{Yuqian2015} investigated which stability solution concepts maximize the probability that the allocations will be stable in an \emph{ a posteriori} fashion, that is after the samples of the uncertainty have been revealed. Finally, \cite{Repeated2} and \cite{Repeated1} addressed uncertainty by studying the dynamics of repeated stochastic
coalitional games. \par
 Unlike the aforementioned works, we construct a data-driven methodology that allows the provision of distribution-free guarantees on the stability of allocations in a probably approximately correct (PAC) fashion, i.e., ensure that, with high confidence, the agents' probability to defect from the grand coalition and form subcoalitions is bounded by a prespecified threshold. Specifically, this bound depends on the amount of available data, the confidence parameter and which samples from the data are crucial to obtain this particular allocation or set of allocations. Connecting PAC learning with uncertain coalitional games has also been considered in \cite{Balcan2015}. Therein, the authors use a sampling-based approach to learn the value function using a randomly generated subset out of the total number of potential coalitions. 
  The spirit of \cite{Procaccia} is similar, using Vapnik-Chervonenkis (VC) theory to learn the winning coalitions for the class of the so-called \emph{simple} games.  \par
Our main contributions in comparison with the existing literature are as follows. Focusing on coalitional games with uncertain value functions, we leverage recent results from the scenario approach  \cite{Ramponi2018,campi2020scenario,campi2021theory} to provide distribution-free PAC-type stability guarantees for agents'  allocations. 
  With respect to \cite{Balcan2015}, our results do not suffer from the conservatism associated with VC-theoretic results. Furthermore, \cite{Balcan2015} focuses on a complementary problem where only a randomized subset of possible coalitions is considered; in contrast, we include all of them but account in a randomized manner for uncertainty in the value functions. \par 
 Our analysis is based on mild assumptions, as we assume no prior knowledge of the sample space or the probability distribution of the uncertainty. Note that using robust versions of the core that take into account all possible uncertainty realizations (see for example \cite{Raja}), as our main solution concept for stability would pose significant challenges in such a general setting. To start with, the evaluation of the robust core under our set-up constitutes a challenging task as the uncertainty support set $\Delta$ is in general unknown in the definition of the robust core (Definition \ref{stable_allocations}); even in cases where this is known, finding the maximum coalitional value could be computationally unviable  when $\Delta$ is a continuous set, unless we impose additional structural assumptions on its geometry.  Furthermore, the robust versions of the core do not provide flexibility, as \cite{Raja} considers the highest value realizations that can possibly occur for each coalition. This unveils the conservatism that accompanies such worst-case paradigms, as it is not possible to account for cases where allocations are allowed to be unstable with a small probability as common in a more general chance-constrained framework. These challenges are circumvented in our work by adopting a data-driven paradigm with the scenario approach as its backbone. 
Our first contribution hence involves the introduction of a data-driven stability concept, i.e., the \emph{scenario core}. Our approach allows then to provide stability guarantees collectively for all allocations that belong to this randomized solution set. \par
A significant impediment on allocation stability that spans the literature of cooperative game theory is the fact that the core of a game, i.e., the set of stable allocations, can be empty. In such cases, a relaxed deterministic version of the core is usually proposed to provide stable solutions by charging the formation of alternative sub-coalitions.  This set is known in the literature as the \emph{$\epsilon$-core} of the game, where $\epsilon$ is a parameter that determines the level of relaxation. In certain cases, the $\epsilon$-core can be viewed as a set of nearly stable solutions. Our work generalizes this concept in a data-driven framework.  Leveraging recent results on the scenario approach with relaxation \cite{campi2020scenario,campi2021theory}, we formulate the problem of finding an allocation in a randomized $\epsilon$-core of a game and we provide probabilistic stability guarantees for this allocation. As a byproduct of the proposed approach, the emptiness of the original core can be easily revealed along with the part of the data that is responsible for the emptiness of the scenario core. For notational reasons, we will refer to the relaxed core as the $\zeta$-core in the remainder of the document. \par  
The rest of the paper is organized as follows. Section II formulates the problem under study establishing a data-driven framework for uncertain coalitional games. Section III provides then probabilistic stability guarantees for allocations that belong to the scenario core and, in case the scenario core is empty, for those allocations inside the randomized $\zeta$-core. Section IV corroborates our theoretical findings by means of a numerical example, while Section V concludes the paper and proposes future research directions of potential interest.
\section{A data-driven approach to uncertain coalitional games} \label{formulation}
\subsection{Problem formulation and stability of allocations}
Let $\pazocal{N}=\{1, \dots, N\}$ be the index set of $N$ agents and consider any subset $S \subseteq \pazocal{N}$ as a \emph{coalition}. In the proposed setting, selfish agents have the incentive to form coalitions in order to achieve a higher individual payoff, or as the means to perform a certain task. The value of a coalition $S \subseteq \pazocal{N}$ is represented by the so-called \emph{value function}, denoted by $u(S)$, $S \subseteq \pazocal{N}$. 
 Exogenous uncertainties affecting a system of agents are prominent in applications \cite{Raja, Fele_coalitional, Malcolm2, Mccain2008}, e.g, in  economic systems, collaborations of companies can be affected by sociopolitical factors, while in energy markets electricity price fluctuations can affect the coalitional values related to the collective minimization of electricity cost of peer-to-peer coalitions.\par  
We consider in this paper a model of a coalitional game where the underlying network is affected by some exogenous uncertainty $\delta$ that takes values in a set $\Delta$ according to a probability distribution $\mathbb{P}$. Note that our model considers $\Delta$ and $\mathbb{P}$ to be fixed but possibly unknown, thus managing to capture several external factors that could potentially lead to changes of the coalitions' values. In this setting, the value functions are extended to their uncertain counterpart $u: 2^\pazocal{N} \times \Delta \rightarrow \mathbb{R}$ which,  given a coalition $S \subset \pazocal{N}$ and an uncertainty realization $\delta \in \Delta$, returns a scalar value. The value function realizations for any subcoalition could be either available from historical data or extracted from some synthetic dataset obtained by some prediction model. In our setting we consider the value function of the grand coalition to be deterministic.  
An \emph{ uncertain game} is then defined as the tuple $G_\Delta=(\pazocal{N}, u, \Delta, \mathbb{P})$.  
 A vector  $\bs{x}=\col((x_i)_{i \in \pazocal{N}}) \in \mathbb{R}^N$, where $x_i$ is the payoff received by agent $i$ is called an \emph{allocation}.  An allocation is \emph{strictly feasible for a coalition} $S$ if $\sum_{i \in S}x_i < u(S)$. Feasibility of a coalition implies that agents have an incentive to form this coalition. The so-called grand coalition $\pazocal{N}$ is called \emph{efficient} if $\sum_{i \in \pazocal{N}}x_i = u(\pazocal{N})$. 
 In this work we are interested in finding efficient allocations (in the grand coalition sense) that are not strictly feasible by any other coalition. Such allocations are called \emph{stable}, as there are no incentives for agents to form coalitions different from the grand one. The set of all stable allocations is called the \emph{core} of the game.
 However, in our setting, the value function is considered to be uncertain, and hence, the standard definition of the core \cite{Chalkiadakis} is not sufficient to capture the desired stability properties. To this end, we extend the notion of the core to account for the presence of uncertainty, as in the following definition.

 \begin{definition}\textup{(Robust core)} \label{stable_allocations} 
 	The \emph{robust core} $C(G_\Delta)$ of an uncertain game $G_\Delta$ is given by $C(G_\Delta)=\{\bs{x} \in  \R^N  : \textstyle\sum_ {i \in \pazocal{N}}x_i = u(\pazocal{N}),   \ 
 		\textstyle\sum_{i \in S}x_i \geq \max_{\delta \in \Delta} u(S, \delta) \text{ for all }  S \subset \pazocal{N} \nonumber \}.$ \hfill$\square$
 \end{definition}

Definition \ref{stable_allocations} shares a similar spirit with Definition 6 in \cite{Raja}. The robust core $C(G_\Delta)$ provides the coalitional game under study with a measure of robust stability in the sense that, for any allocation $\bs{x} \in C(G_\Delta)$, agents have no incentive to defect from the grand coalition to originate sub-coalitions. \par
\subsection{A PAC-learning approach to allocation stability} 
  Unfortunately, computing explicitly the robust core is hard, as we assume no knowledge on the uncertainty support $\Delta$ and the underlying probability distribution $\mathbb{P}$.
 To circumvent this challenge, we adopt a data-driven methodology and approximate the robust core by drawing a finite number of $K$ independent and identically distributed  (i.i.d.) samples $\delta_K \coloneqq (\delta^{(1)}, \ldots, \delta^{(K)} ) \in \Delta^K$, where $\Delta^K$ denotes the cartesian product consisting of $K$ copies of $\Delta$. We refer to vectors $\delta_K$ as multi-samples. This gives rise to the scenario  game $G_K=(\pazocal{N}, u, \delta_K )$. The core of $G_K$, referred to as the \emph{scenario core}, is then defined as  $C(G_K)=\{\bs{x} \in \R^N : \textstyle\sum_ {i \in \pazocal{N}}x_i= u(\pazocal{N}), \nonumber \\ \textstyle\sum_{i \in S}x_i \geq \max\limits_{k=1,\dots, K} u(S, \delta^{(k)}) \ \text{for all} \ S \subset \pazocal{N} \}.$
Furthermore, we impose the following assumption:
\begin{assumption} \label{nonempty}
	For any multi-sample $\delta_K \in \Delta^K$, the scenario core $C(G_K)$ is non-empty. \hfill$\square$
\end{assumption}

 Assumption \ref{nonempty} implies the existence of stable allocations. We will investigate how to waive this assumption for cases where the presence of uncertainty results in an empty scenario core in Section \ref{PAC-section}. On the basis of available data, we then wish to provide guarantees on the probability that payoff allocations $\bs{x} \in \mathbb{R}^N$ inside the scenario core will remain stable even for future, yet unseen, uncertainty realizations. Borrowing tools from the scenario approach  \cite{Pantazis2020,Fabiani2020b}, we define two probabilities of instability. The first one denotes the probability that a certain allocation will become unstable for a new unseen realization of the uncertainty. The second extends this notion to the probability of instability of the entire scenario core.
\begin{definition} \label{violation} \textup{(Probability of instability)}
	\begin{enumerate}
\item Let  $V: \R^{N} \rightarrow [0,1]$. For any $\bs{x} \in \mathbb{R}^N$, we call 
	\begin{align}
		V(\bs{x})= \mathbb{P}\{\delta \in \Delta :  \textstyle\sum_{i \in S}x_i <  u(S, \delta)	\text{ for some } S \subset \pazocal{N}\} \nonumber 
	\end{align}
\emph{probability of allocation instability}.
		\item 	Let $\mathbb{V}: 2^{\R^{N}} \rightarrow [0,1]$. We call 
	 	\begin{align}
		\mathbb{V}(C(G_K))= \mathbb{P}\{\delta \in \Delta :  \exists \ \bs{x} \in C(G_K) :  \textstyle\sum_{i \in S}x_i <  u(S, \delta), \nonumber \\
		 \text{ for some } S \subset \pazocal{N}\} \nonumber 
	   \end{align}
		 \emph{probability of core instability}.   \hfill$\square$
		\end{enumerate}
\end{definition}

The probability of core instability thus denotes the probability that there exist $\delta$ and $S$ with value function $u(S, \delta)$, parameterized by the uncertainty $\delta$, such that at least one of the stable allocations in the scenario core will become unstable, i.e., the agents will defect from the grand coalition to form $S$.  \par
 By leveraging available data, we then aim at bounding with high confidence the probability of core instability, i.e., in a PAC fashion. To achieve this we introduce two key concepts from statistical learning theory, namely the \emph{algorithm} and the \emph{compression set} \cite{MargellosEtAl2015Compression}.
\begin{definition}\textup{(Algorithm)} \label{algorithm} 
	A mapping $A: \Delta^K \rightarrow 2^{\mathbb{R}^N}$ that takes as input a multi-sample  $\delta_K \in \Delta^K$ and returns the scenario core of game $C(G_K)$ is called an \emph{algorithm}. \hfill$\square$
\end{definition} 

 Note that some samples are more important than others in the decision making procedure. In fact, only a subset of $\delta_K$ may be sufficient to produce the same scenario core. As we will see later in the paper, this subset, known as a compression set, dictates the quality of the probabilistic stability guarantees that we can provide and is defined as follows. 
\begin{definition}\textup{(Compression set)} \label{support} 
	With $\mathbb{P}^K$-probability one with respect to the choice of $\delta_K$, a subset  $I \subseteq \{\delta^{(1)}, \dots, \delta^{(K)}\}$  is a \emph{compression set} of $A$ if  $A(\delta_I) = A(\delta_K)$, where $\delta_I$ is a vector whose elements are the samples included in $I$.\hfill$\square$
\end{definition} \par
A compression set of $A$ with minimal cardinality is hence called \emph{minimal} compression set.  A procedure that enumerates such a set is called a  \emph{compression function}.
\section{Probabilistic stability guarantees of allocations} \label{PAC-section}
\subsection{Stability guarantees for a non-empty scenario core} \label{PAC-subsection1}
The following theorem  provides collective guarantees on the stability of allocations in the scenario core against yet unseen value function realizations.
\begin{thm} \label{postability} \textup{(\emph{A posteriori} collective stability guarantees)}
Consider Assumption \ref{nonempty} and algorithm $A$ as in Definition \ref{algorithm} along with its compression function. Fix a confidence parameter $\beta \in (0,1)$ and define the violation level  $\epsilon: \{0,...,K\} \rightarrow [0,1]$ as a function such that
	\begin{align}
		\epsilon(K)=1 \ \text{and} \  \sum_{s=0}^{K-1} {K \choose s}(1-\epsilon(s))^{K-s}=\beta. \label{epsilon}
	\end{align} 
	We then have that 
	\begin{align}
		\mathbb{P}^{K}  \left\{\delta_K \in \Delta^{K} : \mathbb{V}(C(G_K)) \leq  \epsilon(s_K)  \right\} \geq 1-\beta \label{guarantees}
	\end{align}
	where $\mathbb{P}^{K}= \prod_{k=1}^K \mathbb{P}$  and $s_K$ is the cardinality of the minimal compression set.  \hfill$\square$
\end{thm}
\begin{proof}
The proof follows by re-adapting \cite[Th.~6]{fele-a}.
\end{proof}
Roughly speaking, Theorem \ref{postability} states that, with confidence at least $1-\beta$, the probability that a new, yet unseen, uncertainty sample will make an allocation in the scenario core unstable is bounded by $\epsilon$, a function of $s_K$.
	\begin{figure}
		\centering
		\includegraphics[height=4.7cm,width=8.3cm]{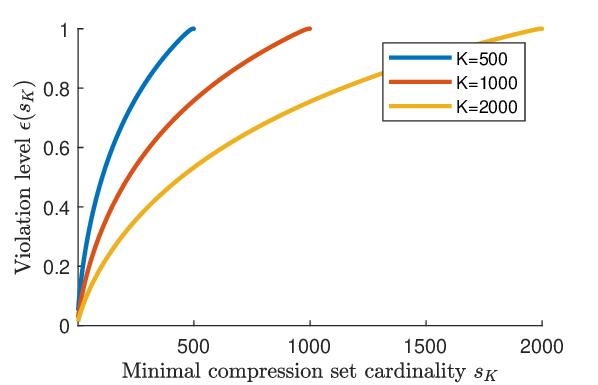}
		\caption{The violation level $\epsilon(s_K)$ as a function of the cardinality of the minimal compression set for $\beta=10^{-6}$ and three different multi-sample sizes, i.e., $K=500, 1000, 2000$. }
		\label{fig1}
	\end{figure}
	 A simple choice of $\epsilon(s)$ can be obtained by splitting $\beta$ evenly among
the $K$ terms in the sum of \eqref{guarantees}, and solving with respect to $\epsilon(s)$. This results in \cite[Eq.~(7)]{Ramponi2018}, and is illustrated numerically in Fig.~\ref{fig1}.
 Note that the nature of Theorem \ref{postability} is \emph{a posteriori}, i.e., we can claim the probabilistic bound in (\ref{guarantees}) only once the samples have been revealed. In fact, this  is required to quantify both the solution and the cardinality $s_K$ of the minimal compression set. \par

  Algorithm~\ref{compression}, introduced in the next page, plays the role of the compression function used to compute $s_K$. Note that the problem at hand has a similar structure to the problem in \cite{Fabiani2020}. Algorithm~\ref{compression}  involves solving a linear optimization program (specifically, a feasibility problem) for each coalition $S' \subset \pazocal{N}$, where we enforce in step 3 the constraint $\sum_{i \in S'} x_i = \max_{k=1, \dots, K} u(S', \delta^{(k)})$, while all other constraints remain unchanged. 
If the problem is feasible, the sample that maximized the right-hand side of this equality constraint is included in the compression set $I$, as feasibility in this case implies that this sample belongs to the compression set of the algorithm $A$ that forms the scenario core $C(G_K)$. Using the union in step 5 ensures that no sample is counted more than once when we compute the cardinality $s_K$ of the compression set $I$ in step 8. 
 \par 
Instead of following the \emph{a posteriori} methodology of Theorem \ref{postability}, one can  provide collective guarantees on the stability of the scenario core also in an \emph{a priori} fashion, whenever the number of samples is greater than the number of possible subcoalitions. As shown in Lemma \ref{apriori_bound}, the \emph{a priori} bound on the cardinality of the compression set coincides with the number of possible coalitions that can be formed; if constraints on the latter are in place, tighter bounds could be provided. The following lemma establishes this observation.
\begin{lemma} \label{apriori_bound}
Assume that $K > |2^{\pazocal{N}}|$.
	The cardinality of the minimal compression set is bounded by the number of possible subcoalitions.
	\hfill$\square$
	\end{lemma}
\emph{Proof:} By the definition of $C(G_K)$, for each $S \subset \pazocal{N}$, the largest value on the right-hand side of the inequality $\sum_{i \in S}x_i \geq  \max_{k=1,\dots, K} u(S, \delta^{(k)})$ is attained at a given sample (assuming no degenerate instances where multiple maximizers exist). Assuming $K > |2^{\pazocal{N}}|$, in a worst-case setting, for each $S \subset \pazocal{N}$ the maximizing sample may be different.
Therefore, for each $S \subset \pazocal{N}$ each inequality $\sum_{i \in S}x_i \geq  \max_{k=1,\dots, K} u(S, \delta^{(k)})$ is a randomized facet of the non-redundant polytope $C(G_K)$. Each of these facets corresponds to exactly one coalition, irrespective of the number of samples. This observation concludes the proof. \hfill $\blacksquare$

Under Lemma \ref{apriori_bound}, by \cite[Theorem 2]{MargellosEtAl2015Compression}, for a fixed confidence level $\beta \in (0,1)$, we can claim that
\[
\mathbb{P}^{K}  \left\{\delta_K \in \Delta^{K} : \mathbb{V}(C(G_K))> \epsilon \right\} \leq \beta,
\]
where $\epsilon$ can be chosen \emph{a priori} such that ${K \choose d} (1-\epsilon)^{K-d} = \beta$, with $d$ being the number of possible subcoalitions. However, note that for games with a high number of agents, the \emph{a posteriori} result of Theorem \ref{postability} might be preferable as the \emph{a priori} one tends to be conservative.

\begin{algorithm}[h]
	\caption{Compression function of algorithm $A$} \label{compression}
	\begin{algorithmic}[1]
		\State	\textbf{initialization:} $I=\emptyset $
		\State \textbf{input:} $\delta_K$
		\State \textbf{for all}  $S' \subset \pazocal{N}$ 
		\begin{equation}	
			\bs{x}^*_{S'} = \left\{
			\begin{aligned}
				& \hspace{0.4cm}  \textrm{arg}\underset{\bs{x} \geq 0}{\textrm{min }} && 0 \\ 
				& \textrm{subject to} &&\sum_{i \in \pazocal{N}}x_i=u(\pazocal{N}), \\
				&&& \sum_{i \in S'} x_i = \max_{k=1, \dots, K} u(S', \delta^{(k)}), \\  
				&&&\sum_{i \in S} x_i \geq \max_{k=1, \dots, K}u(S, \delta^{(k)}) \ \forall \ S \subset \pazocal{N}.  
			\end{aligned}
			\right. \label{feas_coop}
		\end{equation}
		\State \textbf{if} $\bs{x}^*_{S'} \neq \emptyset$%
		\State \hspace{2cm}\ $I \leftarrow I \cup \textrm{arg} \underset{k=1, \dots, K}{\textrm{max }} u(S', \delta^{(k)})$
		\State \textbf{endif}
		\State  \textbf{endfor}
		\State $s_K=|I|$
		\State \textbf{output:} $s_K$
	\end{algorithmic}
\end{algorithm}
\subsection{The case of empty core} \label{empty_core} 

In our framework there might exist some realization $\delta \in \Delta$ resulting in $C(G_\Delta) = \emptyset$ (and hence the same can happen to its scenario counterpart, $C(G_K)$). Establishing whether $C(G_\Delta)$ is nonempty amounts to solving the feasibility program
$$
\left\{
	\begin{aligned}
		& \hspace{0.8cm} \underset{\bs{x} \geq 0}{\textrm{min }} && 0 \nonumber \\ 
		& \textrm{subject to} && \sum_{i \in \pazocal{N}}x_i=u(\pazocal{N}), \nonumber \\  
		&&&\sum_{i \in S} x_i \geq \max_{\delta \in \Delta}u(S, \delta) \ \text{for all} \ S \subset \pazocal{N}. \nonumber  
	\end{aligned}
\right.
$$To circumvent the possibly infinite cardinality of the  constraint set which must hold for all $\delta \in \Delta$, we adopt a data-driven formulation, thus obtaining the following scenario program
$$
\left\{
\begin{aligned}
	    & \hspace{0.8cm} \underset{\bs{x} \geq 0}{\textrm{min }} && 0 \nonumber \\
		& \textrm{subject to}  && \sum_{i \in \pazocal{N}}x_i=u(\pazocal{N}), \nonumber \\  
		&&& \sum_{i \in S} x_i \geq \max_{k=1, \dots, K}u(S, \delta^{(k)}) \ \text{for all} \ S \subset \pazocal{N}.  \nonumber
	\end{aligned}
\right.
$$
In this case, we lift Assumption \ref{nonempty} on non-emptiness of the scenario core and provide through Theorem \ref{relaxation} guarantees on  the probabilistic stability of an allocation in a relaxed version of the scenario core, the so called \emph{scenario  $\zeta$-core}:
\begin{definition}\textup{(Scenario  $\zeta$-core)} \label{epsilon_core}
Given some $\zeta \geq 0$, the scenario $\zeta$-core of a scenario game $G_K$ coincides with the set $C_\zeta(G_K)=\{\bs{x} \in \R^N : \textstyle\sum_{i \in \pazocal{N}}x_i =  u(\pazocal{N}),	\textstyle\sum_{i \in S}x_i \geq  \max_{k=1, \dots, K} u(S, \delta^{(k)})- \zeta  \text{ for any }  S \subseteq \pazocal{N}  \}. $  \hfill $\square$ 
\end{definition}
Note that the standard notion of the scenario core is recovered as a special case of the  scenario $\zeta$-core by setting $\zeta=0$. In other words, the scenario $\zeta$-core is a set of allocations based on available data, where no agent can improve its payoff by leaving the grand coalition. If it happens, then it must pay a penalty of $\zeta$ for leaving. We first aim at finding a solution inside the \emph{least-core}, i.e., the set where allocations are stable with the minimum value of $\zeta$.  Finding such a solution in our set-up amounts to solving the following tractable  optimization program, formulated  in the spirit of the scenario approach with relaxation \cite{campi2020scenario,campi2021theory},
\begin{equation}\label{eq:relaxed_1}
\left\{
	\begin{aligned}
		& \hspace{0.3cm} \underset{\bs{x} \geq 0,~ \xi \geq 0}{\textrm{min}} &&\sum_{k = 1}^K \xi_k  \\ 
		& \textrm{subject to} &&\sum_{i \in \pazocal{N}}x_i=u(\pazocal{N}),   \\
		&&&\sum_{i \in S} x_i \geq \max_{k=1, \dots, K}u(S, \delta^{(k)})- \xi_k, \\
		&&&\hspace{2cm} \text{for all} \ S \subset \pazocal{N},  k=1, \dots, K. 	
	\end{aligned}
\right.
\end{equation}
Note that the proposed results hold independently of the way (\ref{feas_coop}) and (\ref{eq:relaxed_1}) are solved. In this regard, we point out that they amount to linear programs, hence amenable to distributed computation \cite{Bertsekas1997}, \cite{BoydADMM}. 
Let $\xi^\ast_K := \col(\{\xi^\ast_k\}_{k=1}^K)$ and $\zeta^*=\max_{k=1, \dots, K}\xi^*_k$. A solution to (4) is a pair $(\bs{x}^*_K,\xi^\ast_K)$, where $\bs{x}^*_K$ is an allocation in the $\zeta^*$-core (the least core). We make use of the following assumptions: 
 \begin{assumption} \label{uniqueness} \textup{(Uniqueness \cite{campi2020scenario})}
 	For any multi-sample $\delta_K \in \Delta^K$, the solution $(\bs{x}^\ast_K, \xi^\ast_K)$ of \eqref{eq:relaxed_1} is unique. \hfill$\square$
 \end{assumption}
In case the solution is not unique, a solution can be singled out by applying a convex tie-break rule \cite{campi2020scenario}.
 \begin{assumption} \label{nonaccumulation} \textup{(Non-accumulation \cite{campi2020scenario})}
 	For every allocation $\bs{x} \in \mathbb{R}^N$, $\mathbb{P}\{\delta \in \Delta : \sum_{i \in S}x_i =  u(S, \delta) \}=0$. \hfill$\square$
 \end{assumption}

Assumption \ref{nonaccumulation} is related to nondegeneracy \cite{Campi_WaitandJudge} and it is often satisfied when the uncertain parameter $\delta$ does not accumulate, e.g., when $\delta$ admits a density function \cite{campi2020scenario}. From \cite[Thm. 1]{campi2021theory}, consider the polynomial equation
	\begin{align}
		&\binom{K}{s}t^{K-s}-\frac{\beta}{2N}\sum\limits_{i=s}^{K-1}\binom{i}{s}t^{i-s}-\frac{\beta}{6K}\sum\limits_{i=K+1}^{4K}\binom{i}{s}t^{i-s}=0, \label{eq1}\\
		&\text{for any $s=0, \dots, N-1$, while for $s=K$ consider} \nonumber \\
		&1-\frac{\beta}{6K}\sum\limits_{i=K+1}^{4K}\binom{i}{K}t^{i-K}=0. \label{eq2}
	\end{align}
 For any $s=0, \dots, N-1$, equation (\ref{eq1}) has exactly two solutions in $[0, +\infty)$ (see \cite[Thm. 1]{Garatti_RiskComplexity}), which we denote as $\underline{t}(s)$ and $\overline{t}(s)$ with $\underline{t}(s) \leq \overline{t}(s)$. Equation (\ref{eq2}) has only one solution in $[0, +\infty)$ which we denote as $\overline{t}(K)$, while we define $\underline{t}(K)=0$. Let $\underline{\epsilon}(s)=\max\{0, 1-\overline{t}(s)\}$ and $\overline{\epsilon}(s)=1-\underline{t}(s)$, $s=0,1, \dots, K$.
We are now ready to introduce our main theorem on the provision of stability guarantees for allocations inside a scenario $\zeta$-core.
 \begin{thm} \label{relaxation}
 	Consider Assumptions \ref{uniqueness} and \ref{nonaccumulation} and fix a confidence parameter $\beta \in (0,1)$.
 	We then have that 
 	\begin{align}
 		\mathbb{P}^K\{\delta_K \in \Delta^K : \underline{\epsilon}(s^\ast) \leq V(\bs{x}^*_K) \leq \overline{\epsilon}(s^\ast)\} \geq 1- \beta,
 	\end{align}
 	where $\bs{x}_K^*$ is obtained from \eqref{eq:relaxed_1},  $s^\ast$ is the number of $\delta^{(k)}$'s for which there exists  a coalition $S \subset \pazocal{N}$ such that  $\sum_{i \in S} x^*_i \leq u(S, \delta^{(k)})$, and $\underline{\epsilon}(\cdot)$, $\overline{\epsilon}(\cdot)$ are as defined below \eqref{eq2}. \hfill$\square$
 \end{thm}
 \begin{proof}
Note that  for any $\delta \in \Delta$  we can equivalently rewrite $\sum_{i \in S} x_i \geq u(S, \delta)$ for all $S \subset \pazocal{N}$ as $\textrm{max}_{\ell = 1, \ldots, 2^N-1} \, (b_\ell(\delta) - A_\ell \bs{x}) \leq 0$, where the $\ell$-th row $A_\ell$ of the matrix $A \in \R^{(2^N-1) \times N}$ selects those agents generating the $\ell$-th sub-coalition $S_\ell$ among the possible $2^N-1$, while $b(\delta) \coloneqq [u(S_1, \delta) \;  \cdots \; u(S_{2^N-1}, \delta)]^\top \in \R^{2^N-1}$. As such, the program in \eqref{eq:relaxed_1} takes the form
\begin{equation}\label{eq:relaxed_2}
\left\{
\begin{aligned}
 		&\underset{\bs{x} \geq 0, \xi \geq 0}{\textrm{min }} &&\sum_{k = 1}^K \xi_k \\ 
 		&\textrm{ s.t. } &&\sum_{i \in \pazocal{N}}x_i=u(\pazocal{N})   \\
 		&&& f(\bs{x}, \delta^{(k)}) \leq \xi_k, \text{ for all } k = 1,\ldots,K,
	\end{aligned}
\right.
\end{equation}
with $f(\bs{x}, \delta) \coloneqq \textrm{max}_{\ell = 1, \ldots, 2^N-1} \, (b_\ell(\delta) - A_\ell \bs{x})$.
Since, in our setting, \cite[Ass.~1]{campi2021theory} is satisfied and considering  Assumptions \ref{uniqueness} and \ref{nonaccumulation}, direct application of \cite[Th.~1]{campi2021theory} to \eqref{eq:relaxed_2} allows the provision of probabilistic stability guarantees for an allocation $\bs{x}^\ast_K$ inside the scenario least core obtained once the values of $\{\xi^\ast_k\}_{k=1}^K$ have been computed. 
\end{proof}
The formulation in \eqref{eq:relaxed_1} brings several benefits. Specifically, after solving \eqref{eq:relaxed_1}, which is always feasible, if each $\xi^\ast_k = 0$, then the scenario core $C(G_K)$ is nonempty, meaning that, for the collected data, the grand coalition admits at least a stable solution whose robustness can be quantified following the discussion in Subsection \ref{PAC-subsection1} or the bound in Theorem \ref{relaxation}. However, no conclusion can be drawn on the nonemptiness of the original robust core, $C(G_\Delta)$. In case there exists some $\xi^\ast_k > 0$, this means that the scenario core $C(G_K)$ is empty and the same holds for the robust core $C(G_\Delta)$ (since $C(G_\Delta) \subseteq C(G_K)$). Interestingly, Theorem \ref{relaxation} allows us to quantify the penalty to be imposed on the formation of sub-coalitions in order to obtain a solution with a provable probabilistic bound on stability. 
\section{Numerical Study}
\begin{figure}
 	\centering
 	\includegraphics[height=5.2cm,width=8.2cm]{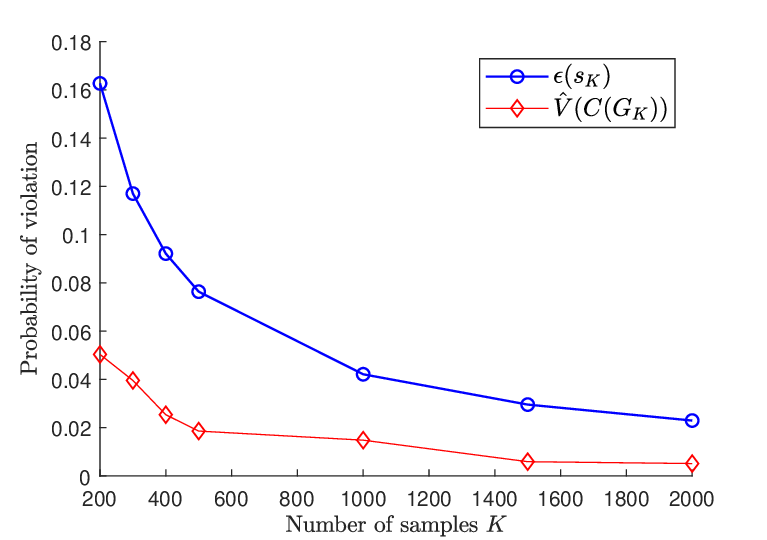}
 	\caption{The violation level $\epsilon(s_K)$ (blue line) vs the empirical probability of violation  $\hat{V}(C(G_K))$ (red line)  for the scenario core of the uncertain game as defined by Table I.  }
	\label{fig2}
 \end{figure}

We consider an uncertain coalitional game with $N=4$ agents and the value functions for each coalition summarized in Table \ref{table}. 
\begin{table} 
	\caption{Uncertain value functions for each coalition} \label{table}
	\centering
	\resizebox{7cm}{!}{
	\begin{tabularx}{0.47\textwidth} { 
			| >{\raggedright\arraybackslash}X 
			| >{\centering\arraybackslash}X 
			| >{\raggedleft\arraybackslash}X | }
		\hline
	\textbf{	Coalition}  &  \textbf{ Value of coalition}   \\
		\hline
		\hline
		$\{1\}$  & $1+ \delta_1$  \\
		\hline
		$\{2\}$ & $1.5+\delta_2$  \\
		\hline
		$\{3\}$  & $1+\delta_3$    \\
		\hline
		$\{4\}$ & $2+\delta_4$ \\
		\hline
		$\{1,2\}$  & $6.5+\delta_{12}$    \\
		\hline
		$\{2,3\}$  & $6.5+\delta_{23}$    \\
		\hline
		$\{3,4\}$ & $7+\delta_{34}$ \\
		\hline
		$\{1,3\}$  & $6+\delta_{13}$    \\
		\hline
		$\{1,4\}$  & $7+\delta_{14}$    \\
		\hline
		$\{2,4\}$ & $7.5+\delta_{24}$ \\
		\hline
		$\{1,2,3\}$  & $11.5+\delta_{123}$    \\
		\hline
		$\{1,2,4\}$  & $12.5+\delta_{124}$    \\
		\hline
		$\{1,3,4\}$ & $12+\delta_{134}$ \\
		\hline
		$\{2,3,4\}$  & $12.5+\delta_{234}$    \\
		\hline
		$\{1,2,3,4\}$ & $17.3$    \\
		\hline 
	\end{tabularx}}
\end{table}
 Each element of the uncertainty $\delta=\{\delta_S\}_{S \subset \pazocal{N}}$ that affects the coalitional values is assumed to be drawn according to a uniform distribution with support in $[-0.5,0.5]$. 
 Drawing a total number of $K= \{200, 300, 400, 500, 1000, 1500, 2000 \}$ samples we use Algorithm \ref{compression} to find the cardinality of the compression, which for this particular example is always upper bounded by $4$. Fixing the confidence parameter to $\beta=10^{-4}$ and considering $10^6$ test samples, we calculate the empirical probability of violation $\hat{V}(C(G_K))$ by checking which of these test samples give rise to values of coalitions that render at least one allocation in the scenario core unstable. Then, we compare in Figure \ref{fig2} this probability with the \emph{a posteriori} theoretical violation level $\epsilon(s_K)$ of the scenario core for each value of $K$. Note that the empirical probability of violation is always below the theoretical level, thus verifying the result of Theorem 1 numerically. \par
\begin{figure}
	\centering
	\includegraphics[height=5.5cm, width=8.5cm]{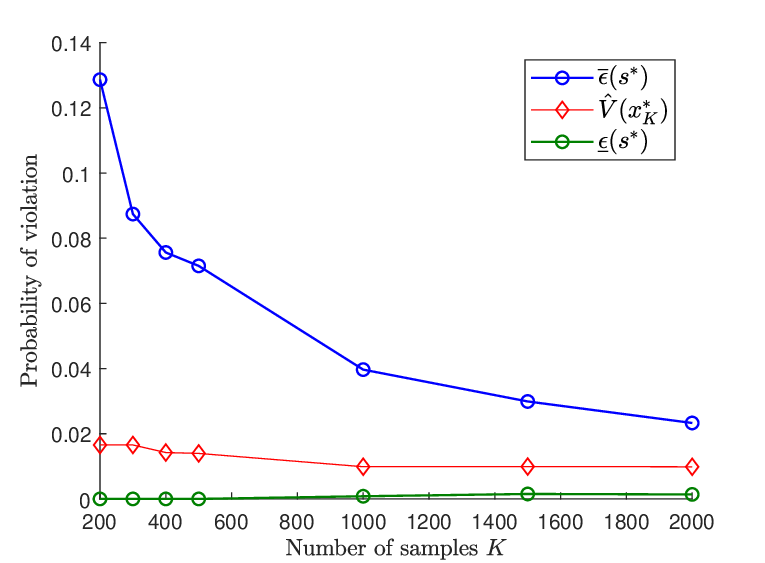}
	\caption{The violation levels $\overline{\epsilon}(s^*), \underline{\epsilon}(s^*)$ vs the empirical probability of violation  $\hat{V}(\bs{x}_K^*)$  for a solution inside the  the scenario $\zeta$-core of the considered uncertain game.}
	\label{fig3}
\end{figure}

Next, we focus on instances of games where the scenario core as defined in Section \ref{formulation} is empty due to the presence of uncertainty. We consider the same coalitional game as before where the uncertainty $\delta$ follows instead a gaussian distribution with mean 0 and standard deviation 0.3 truncated at [-1.5, 1.5]. In this case, the scenario core is empty for a large enough number of samples. Using the relaxation technique of Subsection \ref{empty_core} we find a solution that lies inside the randomized least core. Then, Theorem \ref{relaxation} allows to provide guarantees on the stability of an allocation $\bs{x}^*_K$ in this relaxed core.
 Drawing a different number of samples $K$ each time from the set $\{200, 300, 400, 500, 1000, 1500, 2000\}$ and fixing $\beta=10^{-5}$, we calculate $s^*$ for each multi-sample size and 
	 compute the lower and upper violation levels $\underline{\epsilon}(s^*), \overline{\epsilon}(s^*)$, respectively.
Using a total number of $10^6$ test samples we calculate the empirical probability $\hat{V}(\bs{x}_K^*)$ and we compare it  with the theoretical bounds $\overline{\epsilon}(s^*), \underline{\epsilon}(s^*)$ (these are obtained as discussed below (6), using the code made available in \cite{Garatti_RiskComplexity}) for each multi-sample size. The results are summarized in Figure \ref{fig3}. In general, a lower probability of violation is expected by increasing the number of samples. However, note that a larger number of samples implies in our case that a larger number of relaxed constraints is generally needed to find an allocation by solving (\ref{eq:relaxed_1}). Thus, such allocations in the relaxed core have similar or only slightly decreasing empirical probabilities of violation as the number of samples increases. 
 \section{Conclusion}
 Focusing on coalitional games with uncertain value functions, we borrow tools from the scenario theory to provide distribution-free guarantees on the probability that agents will not defect from the grand coalition to originate subcoalitions for unseen uncertainty realizations. We provided collective stability guarantees for the entire core set in a PAC learning fashion, and accounted for the case where the scenario core can be empty, proposing a methodology to accompany allocations in a relaxed core with stability certificates. Future work will focus on establishing tighter stability guarantees on a specific core subset where desired qualities of the solution, e.g., fairness within a certain amount of tolerance, are taken into account.
\bibliographystyle{IEEEtran}
\bibliography{biblio_coop}
\end{document}